\documentclass[11pt]{article}
\usepackage{amsthm}
\usepackage[nocompress]{cite}
\usepackage[dvips]{graphics,color}
\usepackage{amssymb}
\usepackage{amsfonts}
\usepackage{graphicx}
\usepackage{subfigure}
\usepackage{hyperref}
\usepackage{algorithm}
\usepackage{algorithmic}
\usepackage{geometry}
\usepackage{epsfig}
\usepackage{listings}

\begin{document}

\title{Harmonic and Refined Harmonic Shift-Invert Residual Arnoldi and
Jacobi--Davidson Methods for Interior Eigenvalue Problems\footnote{Supported
by National Basic Research Program of China 2011CB302400 and the
National Science Foundation of China (No. 11071140).}}

\author{Zhongxiao Jia\thanks{Department of Mathematical Sciences,
Tsinghua University, Beijing 100084, People's Republic of China,
{\sf jiazx@tsinghua.edu.cn}.} \and Cen Li\thanks{Department of Mathematical
Sciences, Tsinghua University, Beijing 100084, People's Republic of China,
{\sf licen07@mails.tsinghua.edu.cn}.}}
\date{}
\maketitle

\newtheorem{Def}{\bf Definition}
\newtheorem{Exm}{\bf Example}
\newtheorem{theorem}{\bf Theorem}
\newtheorem{Lem}{\bf Lemma}
\newtheorem{Cor}{\bf Corollary}
\newtheorem{Alg}{\bf Algorithm}
\newtheorem{Exp}{\bf Experiment}
\newtheorem{Ass}{\bf Assumption}
\newtheorem{Rem}{\bf Remark}
\newtheorem{Pro}{\bf Proposition}

\newcommand{\Span}{\mathrm{span}}
\newcommand{\xp}{\mathbf{x}_\perp}
\newcommand{\Xp}{X_\perp}
\newcommand{\sep}{\mathrm{sep}}
\newcommand{\diag}{\mathrm{diag}}
\newcommand{\imag}{\mathrm{i}}
\newcommand{\bfSig}{\mathbf{\Sigma}}
\newcommand{\mV}{\mathcal{V}}
\newcommand{\wopt}{w_{\mathrm{opt}}}
\newcommand{\Ran}{\mathrm{Ran}}

\begin{abstract}
This paper concerns the harmonic shift-invert residual Arnoldi (HSIRA)
and Jacobi--Davidson (HJD) methods as well as their refined
variants RHSIRA and RHJD for the interior eigenvalue problem. Each method needs
to solve an inner linear system to expand the subspace successively.
When the linear systems are solved only approximately, we are led
to the inexact methods. We prove that the inexact
HSIRA, RHSIRA, HJD and RHJD methods mimic their exact counterparts
well when the inner linear systems are solved with only {\em low} or
{\em modest} accuracy. We show that
(i) the exact HSIRA and HJD expand subspaces better than the exact SIRA and JD
and (ii) the exact RHSIRA and RHJD expand subspaces better than the exact
HSIRA and HJD. Based on the theory, we design stopping criteria for inner solves.
To be practical, we present restarted HSIRA, HJD, RHSIRA and RHJD algorithms.
Numerical results demonstrate that these algorithms are much more
efficient than the restarted standard SIRA and JD algorithms
and furthermore the refined harmonic algorithms outperform the harmonic
ones very substantially.
\smallskip

\smallskip

\end{abstract}

\section{Introduction}

Consider the linear eigenproblem
\begin{equation}
Ax=\lambda x,\ \ \ x^Hx=1,
\end{equation}
where $A\in\mathcal{C}^{n\times n}$ is large and possibly sparse
with the eigenvalues labeled as
\begin{eqnarray*}
0<|\lambda_1-\sigma|<|\lambda_2-\sigma|\leq\cdots\leq|\lambda_n-\sigma|,
\end{eqnarray*}
where $\sigma\in\mathcal{C}$ is a given target inside the spectrum of $A$.
We are interested in the eigenvalue
$\lambda_1$ closest to the target $\sigma$ and/or the associated eigenvector
$x_1$. This is called the interior eigenvalue problem.
We denote $(\lambda_1,x_1)$ by $(\lambda,x)$ for simplicity.

The interior eigenvalue problem arises from many applications \cite{matrixmarket}.
Projection methods have been widely used for it
\cite{bai2000templates,parlett1998symmetric,saad1992eigenvalue,
vandervorst2002eigenvalue,stewart2001eigensystems}.
For a given subspace $\mathcal{V}$, the standard projection
method, i.e., the Rayleigh--Ritz method,  seeks the approximate
eigenpairs, i.e., the Ritz pairs,  $(\nu,y)$ satisfying
\begin{equation}\label{stanpro}
(A-\nu I)y\perp\mathcal{V},\ \ \ y\in\mathcal{V}.
\end{equation}
As is well known, the standard projection method is quite effective
for computing exterior eigenpairs, but is not for computing
interior ones. For interior eigenpairs, one usually faces two difficulties
when using the standard projection method: (i) the Ritz vectors might be
very poor \cite{jia1995,jia2001analysis} even if the subspace that contains
good enough information on the desired eigenvectors;
(ii) even though there are Ritz values which are good approximations to
the desired eigenvalues, it is usually hard to recognize them and
to select the {\em correct} Ritz pairs to
approximate the desired interior eigenpairs \cite{morgan1991,morgan1998}.
Because of possible mismatches, the Ritz pairs may misconverge or converge
irregularly as subspace $\mV$ is expanded or restarted.

To overcome the above second difficulty, one popular approach is to
apply the standard projection method
to the shift-invert matrix $(A-\sigma I)^{-1}$,
which transform the desired eigenvalues of $A$ into the dominant (exterior) ones.
However, if it is very or too costly to factor $A-\sigma I$,
as is often the case in practice, one must resort
iterative solvers to approximately solve the linear systems involving $A-\sigma I$.
Over years, one has focused the inexact shift-invert Arnold (SIA)
type methods where inner linear systems
are solved approximately. It has turned out that inner linear systems
must be solved with very high accuracy when approximate eigenpairs are
of poor accuracy \cite{spence2009ia,simoncini2005ia,simoncini2003ia}.
As a result, it may be very expensive and even impractical
to implement SIA type methods.

As an alternative, one uses the harmonic projection method, i.e., the harmonic
Rayleigh--Ritz method, to seek the approximate eigenpairs, the
harmonic Ritz pairs, $(\nu,y)$
satisfying
\begin{equation}\label{harpro}
(A-(\nu I)y\perp(A-\sigma I)\mathcal{V},\ \ \ y\in\mathcal{V}.
\end{equation}
The harmonic projection method has been used widely
for the interior eigenvalue problem \cite{bai2000templates,saad1992eigenvalue,
vandervorst2002eigenvalue,stewart2001eigensystems}.
The method can be derived by using the standard Rayleigh--Ritz method
on $(A-\sigma I)^{-1}$ with respect to the subspace $(A-\sigma I)\mathcal{V}$
without factoring $A-\sigma I$. The harmonic projection method has the advantage
that $\frac{1}{\mu-\sigma}$ is the Ritz value of $(A-\sigma I)^{-1}$.
This is helpful to select the {\em correct} approximate eigenvalues since
the desired interior eigenvalues near $\sigma$ have been transformed into
the exterior ones that are much easier to match. So the harmonic projection method
overcomes the second difficulty well.

However, the harmonic projection method inherits the first difficulty since
the harmonic Ritz vectors may still be poor, converge irregularly and even
fail to converge; see \cite{jia2005analysis} and also
\cite{stewart2001eigensystems,vandervorst2002eigenvalue} for a systematical
account. To this end, the refined projection method, i.e.,
the refined Rayleigh--Ritz method, has been proposed in
\cite{jia1997refined,jia2001analysis} that cures the possible non-convergence
of the standard and harmonic Ritz vectors. The refined method replaces
the Ritz vectors by the refined Ritz vectors that minimize residuals formed
with the available approximate eigenvalues
over the same projection subspace. The refined projection principle has also been
combined with the harmonic projection method. Particularly,
the harmonic and refined harmonic Arnoldi methods were proposed in
\cite{jia2002refinedharmonic} that are shown to be more effective for computing
interior eigenpairs than the standard and refined Arnoldi methods.
Due to the optimality, refined (harmonic) Ritz vectors are better approximations
than (harmonic) Ritz vectors.
The convergence of the refined projection methods has been established
in a general setting \cite{jia2005analysis,jia2001analysis}; see also
\cite{stewart2001eigensystems,vandervorst2002eigenvalue}.

A necessary condition for the convergence of any kind of projection methods
is that the subspace contains enough accurate approximations to the desired
eigenvectors. Therefore, a main ingredient for the success of a projection method
is to construct a sequence of subspaces that contain increasingly
accurate approximations to the desired eigenvectors. To achieve this goal,
it is naturally hoped that each expansion vector makes contribution to the desired
eigenvectors as much as possible. It turns out that preconditioning plays
a key role in obtaining good expansion vectors \cite{morgan1998,vandervorst1996jd}.
Such kind of methods includes
the Davidson method \cite{davidson1975,crouzeix1994},
the Jacobi--Davidson (JD) method \cite{vandervorst1996jd}
and shift-invert residual Arnoldi (SIRA) method \cite{leestewart07}.

The Davidson method can be seen as a preconditioning Arnoldi (Lanczos) method
and a prototype of the SIRA and JD methods, in which an inner linear system
involving $A-\sigma I$ needs to be solved at each outer iteration.
When a factorization of $A-\sigma I$
is impractical, only iterative solvers are viable for solving the inner
linear systems. This leads to the inexact SIRA and JD.
A fundamental issue on them is how accurately we should solve inner linear systems
in order to make the inexact methods mimic the exact counterparts well, that is,
the inexact and exact ones use comparable outer iterations to achieve the convergence.
For the simplified or single vector JD method without subspace acceleration,
it has been argued that it may not be necessary to solve inner
linear system  with high accuracy and a low or modest accuracy may
suffice \cite{hochstenbachnotay09,vandervorst1996jd}. This issue is more
difficult to handle for the standard JD method.
Let the relative error of approximate solution of
an inner linear system be $\varepsilon$. Lee and Stewart \cite{leestewart07}
have given some analysis of the inexact SIRA under a series
of assumptions. However, their main results appear hard to justify or interpret,
and furthermore no quantitative estimates for $\varepsilon$ have been given
eventually. Using a new analysis approach, the authors \cite{jiali11} have
established a rigorous general convergence theory of the inexact SIRA and JD in a
unified way and derived quantitative and explicit estimates for $\varepsilon$.
The results prove that a low or moderate $\varepsilon$,
say $[10^{-4},10^{-3}]$, is generally enough to make the inexact SIRA and JD
behave very like the exact methods. As a result, the methods are expected to
be more practical than the inexact SIA method.
The authors have confirmed the theory numerically in \cite{jiali11}.

Because of the merits of harmonic and refined harmonic projection
methods for the interior eigenvalue problem, in this paper
we consider the harmonic and refined harmonic variants,
HSIRA, RHSIRA and HJD, RHJD, of the standard SIRA and JD methods. Exploiting
the analysis approach and some results in \cite{jiali11},
we establish the general convergence theory and derive the estimates for
$\varepsilon$ for the four inexact methods.
The results are similar to those for the SIRA
and JD methods in \cite{jiali11}. They prove that, in order to make the inexact
methods behave like their exact counterparts, it is generally enough
to solve all the inner linear systems with only low or modest accuracy.
Furthermore, we show that the exact HSIRA and HJD expand subspaces better
than the exact SIRA and JD and, in turn, the exact RHSIRA and RHJD expand
subspaces better than the exact HSIRA and HJD. These results and the merits
of the harmonic and refined harmonic methods
mean that the harmonic Ritz vectors and refined harmonic
Ritz vectors are better for subspace expansions and restarting than the standard
Ritz vectors. Based on the theory, we design stopping criteria for inner solves.
To be practical, we present restarted HSIRA, RHSIRA, HJD and RHJD algorithms.
We make numerical experiments to confirm our theory and demonstrate that
the restarted HSIRA, HJD and RHSIRA and RHJD algorithms are considerably
more efficient than the restarted standard SIRA and JD algorithms
and furthermore the restarted RHSIRA and RHJD outperform the others very greatly
for the interior eigenvalue problem.

This paper is organized as follows.
In Section~\ref{sec:harsirajd}, we describe the harmonic and refined
harmonic SIRA and JD methods. In Section~\ref{sec:inexact},
we present our results on $\varepsilon$.
In Section~\ref{sec:restart} we describe restarted algorithms.
In Section~\ref{sec:experiments}, we report numerical experiments
to confirm our theory. Finally, we conclude the paper and point out future
work in Section~\ref{sec:conc}.

Some notations to be used are introduced. Throughout the paper, denote
by $\|\cdot\|$ the Euclidean norm of a vector and the spectral norm of a matrix,
by $I$ the identity matrix with the order clear from the context, and
by the superscript $H$ the complex conjugate transpose of a vector or matrix.
We measure the deviation of a nonzero vector $y$ from a subspace $\mV$ by
$$
\sin\angle(\mV,y)=\frac{\|(I-P_\mV)y\|}{\|y\|},
$$
where $P_\mV$ is the orthogonal projector onto $\mV$.

\section{The harmonic SIRA and JD methods and their refined versions}
\label{sec:harsirajd}


Let the columns of $V$ form an orthonormal basis of a given general subspace $\mV$
and define the matrices
\begin{eqnarray}\label{defHG}
H&=&V^H(A-\sigma I)^HV,\\
G&=&V^H(A-\sigma I)^H(A-\sigma I)V.
\end{eqnarray}
Then the standard projection method computes the Ritz pairs $(\nu,y)$ satisfying
\begin{equation}\label{seig}
\left\{\begin{array}{rcl}
H^Hz &=&\mu z ~~{\rm with}~~\|z\|=1,\\
\nu&=&\mu+\sigma,\\
y&=&Vz,
\end{array}
\right.
\end{equation}
and selects $\nu$'s closest to  $\sigma$, which correspond to the smallest $\mu$'s
in magnitude, and the associated $y$'s to approximate
the desired eigenpairs of $A$. Given an initial subspace $\mV$ and the
target $\sigma$, we describe the standard SIRA and JD methods as Algorithm~1 
for our later use.

\begin{algorithm}\label{alg:sirajd}
\caption{The SIRA/Jacobi--Davidson method with the target $\sigma$}
\begin{algorithmic}
\STATE{Given the target $\sigma$ and a user-prescribed convergence tolerance
$tol$, suppose the columns of $V$ form an orthonormal basis of the initial subspace
$\mathcal{V}$.}
\REPEAT
\STATE{\begin{enumerate}
\setlength{\itemsep}{0pt}
\setlength{\parsep}{0pt}
\setlength{\parskip}{0pt}
\item Use the standard projection method to compute the Ritz pair $(\nu,y)$ of $A$
with respect to $\mV$, where $\nu\cong\lambda$ and $\|y\|=1$.
\item Compute the residual $r=Ay-\nu y$.
\item In SIRA, solve the linear system
\begin{equation}\label{lssira}
(A-\sigma I)u_S=r.
\end{equation}
In JD, solve the correction equation
\begin{equation}\label{lsjd}
\label{ls_jd}(I-yy^H)(A-\sigma I)(I-yy^H)u_J=-r
\end{equation}
for $u_J\perp y$.
\item Orthonormalize the $u_S$ or $u_J$ against $V$ to get the expansion vector $v$.
\item Subspace expansion: $V=\left[\begin{array}{cc} V & v\end{array}\right]$.
\end{enumerate}}
\UNTIL{$\|r\|<tol$}
\end{algorithmic}
\end{algorithm}

The harmonic projection method seeks
the harmonic Ritz pairs $(\mu,y)$ satisfying (\ref{harpro}), which amounts
to solving the generalized eigenvalue problem of the pencil $(H,G)$:
\begin{equation}\label{hareig}
\left\{\begin{array}{rcl} Hz&=&\frac{1}{\mu} Gz~~{\rm with}~~\|z\|=1,\\
\nu&=&\mu+\sigma\\\
y&=&Vz,\end{array}\right.
\end{equation}
and selects $\nu$'s closest to $\sigma$, which corresponds to
the smallest $\mu$'s in magnitude, and the associated $y$'s to approximate
the desired eigenpairs of $A$.

Alternatively, we may use the Rayleigh quotient $\rho$ of $A$
with respect to the harmonic Ritz vector $y$ to approximate $\lambda$ as well.
Recall the definition (\ref{defHG}) of $H$ and note from (\ref{hareig})
that $y=Vz$ and $\|z\|=1$. Then
\begin{equation}\label{rho}
\rho=\frac{y^HAy}{y^Hy}
=\frac{z^HV^H(A-\sigma I+\sigma I)Vz}{z^HV^HVz}=z^HH^Hz+\sigma,
\end{equation}
so we can compute $\rho$ very cheaply.
It is proved in \cite{jia2005analysis} that once $\sigma$ is very close to
the desired eigenvalue $\lambda$, then the harmonic projection method may miss
$\lambda$ and the harmonic Ritz value $\mu$ may misconverge
to some other eigenvalue; see \cite[Theorem 4.1]{jia2005analysis}.
However, whenever $y$ converges to the desired
eigenvector $x$, the Rayleigh quotient $\rho$ must converge to $\lambda$
no matter how close $\sigma$ is to $\lambda$;
see \cite[Theorem 4.2]{jia2005analysis}. Therefore, it is better and safer
to use $\rho$, rather than $\mu$, as an approximate eigenvalue.
Another merit is that $\rho$ is optimal in the sense of
$$
\rho=\arg~\min\limits_{\nu\in\mathcal{C}}\|(A-\nu I)y\|.
$$
In the sequel, as was done in the literature, e.g.,
\cite{jia2005analysis,morgan1998}, we will always use
$\rho$ as the approximation to $\lambda$
in the harmonic projection method.

Define $r=(A-\rho I)y$. Then $(A-\rho I)y\perp y$. We replace
$r$, the residual of standard Ritz pair, in (\ref{lssira}) and (\ref{lsjd})
by this new $r$ and the standard Ritz vector $y$ in (\ref{lsjd}) by the harmonic
Ritz vector. We then solve the resulting (\ref{lssira}) and (\ref{lsjd})
with the new $r$ and $y$, respectively. With the solutions, we expand
the subspaces analogous to Steps 4--5 of Algorithm 1. In such a way,
we get the HSIRA and HJD methods.

We should point out that a harmonic JD method had already been proposed
as early as the standard one in \cite{vandervorst1996jd}, where it is suggested
to solve the following correction equation for $u_J'\perp y'$
(suppose that $y'^Hy\neq0$)
\begin{equation}\label{lsjd'}
\left(I-\frac{yy'^H}{y'^Hy}\right)(A-\sigma I)
\left(I-\frac{yy'^H}{y'^Hy}\right)u_J'=-r
\end{equation}
with $r=(A-\mu I)y$ and $y'=(A-\sigma I)y$. It is different from (\ref{lsjd})
in the harmonic JD method described above, so its solution is different
from that of (\ref{lsjd}) too. We prefer (\ref{lsjd}) since there is
a potential danger that the oblique projector involved in (\ref{lsjd'}) may make
it worse conditioned than (\ref{lsjd}) considerably.


In what follows we propose the RHSIRA and RHJD methods.
We first compute the Rayleigh quotient $\rho$ defined by (\ref{rho}) and then
seek a unit length vector $\hat{y}\in\mV$ satisfying the optimality
property
\begin{equation}\label{refinedy}
\left\|(A-\rho I)\hat{y}\right\|
=\min\limits_{w\in\mV,~\|w\|=1}\left\|(A-\rho I)w\right\|
\end{equation}
and use it to approximate $x$. So $\hat{y}$ is the best approximation to $x$
from $\mV$ with respect to $\rho$ and the Euclidean norm.
We call $\hat{y}$ a refined harmonic Ritz vector or more generally a refined
approximate eigenvector since the pair $(\rho,\hat{y})$ does not
satisfy the harmonic projection (\ref{harpro}) any more.
Since $V$ forms an orthonormal basis of $\mathcal{V}$,
(\ref{refinedy}) is equivalent to finding the unit length
$\hat{z}\in\mathcal{C}^{m}$ such that $\hat y=V\hat z$ with $\hat z$ satisfying
\begin{eqnarray*}
\|(A-\rho I)\hat{y}\|
&=&\min\limits_{z\in\mathcal{C}^m,~\|z\|=1}\left\|(A-\rho I)Vz\right\|\\
&=&\left\|(A-\rho I)V\hat{z}\right\|\\
&=&\sigma_{\min}\left((A-\rho I)V\right).
\end{eqnarray*}
We see that $\hat{z}$ is the right singular vector
associated with smallest singular value of the matrix $(A-\rho I)V$,
and $\hat{y}=V\hat{z}$.
So we can get $\hat{z}$ and $\|(A-\rho I)V\hat{z}\|$
by computing the singular value decomposition (SVD) of $(A-\rho I)V$.

Similar to (\ref{rho}), define
\begin{equation}\label{uprho}
\hat\rho=\hat y^HA\hat y=\hat z^HH^H\hat z+\sigma,
\end{equation}
the Rayleigh quotient of $A$ with respect to
the refined harmonic Ritz vector $\hat y$.
Then the new residual $\hat r=(A-\hat\rho)\hat y\perp \hat y$. Replace $r$ and $y$
by $\hat r$ and $\hat y$ in (\ref{lssira}) and (\ref{lsjd}) and perform
Steps 3--5. Then we obtain the RHSIRA and RHJD methods, respectively.
We will use $(\hat\rho,\hat y)$ to approximate $(\lambda,x)$ in the methods
as $\|(A-\hat\rho I)\hat y\|\leq \|(A-\rho I)\hat y\|$.

Some important results have been established in
\cite{jia2004,jia2005analysis} for standard and refined projection methods.
The following two results in \cite{jia2004} are directly adapted to the
harmonic and refined harmonic projection methods:
First, we have $\|(A-\rho I)\hat{y}\|<\|(A-\rho I)y\|$
unless the latter is zero, that is, the pair $(\rho,y)$
is an exact eigenpair of $A$. Second, if $\|(A-\rho I)y\|\neq0$
and there is another harmonic Ritz value close to $\rho$, then
it may occur that
$$
\|(A-\rho I)\hat{y}\|\ll\|(A-\rho I)y\|.
$$
So $\hat{y}$ can be much more accurate than $y$ as an
approximation to the desired eigenvector $x$.

For a general $m$-dimensional subspace $\mV$, two approaches are proposed
in \cite{jia2000,jia2006svd} for computing $\hat z$.
Approach I is to directly form the cross-product matrix
\begin{eqnarray}\label{defS}
S&=&V^H(A-\rho I)^H(A-\rho I)V\\
\label{comS}&=&G+(\overline{\sigma-\rho})H^H+(\sigma-\rho)H+|\sigma-\rho|^2I,
\end{eqnarray}
which is Hermitian semi-positive definite.
The desired $\hat{z}$ is the normalized eigenvector associated with
the smallest eigenvalue of $S$, and $\|(A-\rho I)\hat{y}\|$ is the square root of
the smallest eigenvalue. Noticing that $G$ and $H$ are already available in the
procedure of the harmonic projection, we can form $S$ using at most $12m^2$ flops
by taking the worst case into account that all terms in (\ref{comS})
are complex. So, as a whole, we can compute $\hat{z}$ by the QR algorithm
at cost of $O(m^3)+12m^2$ flops.

Approach II is to first make the thin or compact QR decomposition $(A-\rho I)V=QR$
and then make the SVD of the $m\times m$ triangular matrix
$R$ by $O(m^3)$ flops to get the smallest singular value and the associated
right singular vector, which are $\|(A-\rho I)\hat{y}\|$ and $\hat{z}$,
respectively. If we use the Gram--Schmidt orthogonalization procedure
with iterative refinement to compute the QR decomposition,
then we will get a numerically orthonormal $Q$
within a small multiple of the machine precision, which totally needs
no more at most $4nm^2$ flops generally if $\rho$ is real.
As a whole, $4nm^2+O(m^3)$ flops are needed.

By comparison, we find that Approach I is, computationally, much more effective
than Approach II. It can be justified \cite{jia2006svd} that
Approach I is numerically stable for computing $\hat z$ provided that
there is a considerable gap between the smallest singular value and the
second smallest one of $(A-\rho I)V$.
So in this paper, we use Approach I to compute $\hat{y}$.

Because of different right-hand sides,
it is important to note that expanded subspaces are generally different for
the SIRA, HSIRA and RSIRA methods whatever the linear systems
are solved either exactly or approximately. The same is true for the
JD, HJD and RHJD methods because of not only different right-hand sides
but also different $y$ in the correction equation (\ref{lsjd}).
This is different from shift-invert Arnoldi (SIA) type methods,
i.e., the standard SIA, harmonic SIA, refined SIA and refined harmonic SIA, where
the the updated subspaces are the same once the linear systems
are solved exactly or approximately using the same inner solver
with the same accuracy because the right-hand sides involved are
always the currently newest basis vector at each outer iteration step.
We will come back to this point
in the end of Section~\ref{sec:inexact} and show which method favors
subspace expansion.

\section{Accuracy requirements of inner iterations in HSIRA, HJD, RHSIRA and
RHJD}
\label{sec:inexact}

In this section, we review some important results in \cite{jiali11} and
apply them to establish the convergence theory of the inexact HSIRA, HJD
and RHSIRA and RHJD. We prove that each inexact method
mimics its exact counterpart well provided that all the inner linear
systems are solved with only low or modest accuracy.
We stress that in the presentation $y$ is a general approximation to $x$.
Returning to our methods, $y$ is just the
harmonic Ritz vector $y$ in HSIRA and HJD
and the refined harmonic Ritz vector $\hat{y}$
in RHSIRA and RHJD. Finally, we look into the way that each exact
method expands the subspace and make a simple analysis, showing
that HSIRA, HJD and RHSIRA, RHJD generally expand the subspace more
effectively than SIRA and JD when computing interior eigenpairs,
so that the former ones are more effective than the latter ones.
This advantage conveys to the inexact HSIRA, HJD and RHSIRA, RHJD
when each inexact method mimics its exact counterpart well.

We can write the linear system (\ref{lssira}) and the correction
equation (\ref{lsjd}) as a unified form
\begin{equation}\label{lsunified}
(A-\sigma I)u=\alpha_1y+\alpha_2(A-\sigma I)y,
\end{equation}
which leads to (\ref{lssira}) in the HSIRA or RHSIRA method and (\ref{lsjd}) in
the HJD or RHJD method when $\alpha_1=\sigma-\nu,\ \alpha_2=1$ or
$\alpha_1=\sigma-\rho,\ \alpha_2=1$ and $\alpha_1=-\frac{1}{y^HBy}, \ \alpha_2=1
\mbox{ or }\alpha_1=-\frac{1}{\hat{y}^HB\hat{y}},
\ \alpha_2=1$, respectively. Define $B=(A-\sigma I)^{-1}$.
Then the exact solution $u$ of (\ref{lsunified}) is
\begin{equation}\label{defu}
u=\alpha_1By+\alpha_2y.
\end{equation}
Let $\tilde{u}$ be an approximate solution of (\ref{lsunified})
and its relative error be
\begin{equation}\label{epsilon}
\varepsilon=\frac{\|\tilde{u}-u\|}{\|u\|}.
\end{equation}
Then we can write
\begin{equation}
\tilde{u}=u+\varepsilon\|u\|f,
\end{equation}
where $f$ is the normalized error direction vector. Note that the direction of $f$
depends on inner linear systems solves and accuracy
requirements, but it shows no particular features in general. So $\sin\angle(\mV,f)$
is generally moderate and not near zero.
The vectors $u$ and $\tilde u$ are orthogonalized against $V$
to get the (unnormalized) expansion vectors $(I-P_\mV)u$
and $(I-P_\mV)\tilde{u}$, respectively,
where $P_\mV$ is the orthogonal projector onto $\mV$.
Define the relative error
\begin{equation}\label{tepsilon}
\tilde\varepsilon
=\frac{\|(I-P_\mV)\tilde{u}-(I-P_\mV)u\|}{\|(I-P_\mV)u\|},
\end{equation}
which measures the relative difference between two expansion vectors $(I-P_\mV)u$
and $(I-P_\mV)\tilde{u}$.
The following result \cite[Theorem 3]{jiali11} establishes a compact bound
for $\varepsilon $ in terms of $\tilde\varepsilon$.

\begin{theorem}\label{lem:epsilon}
Define $B=(A-\sigma I)^{-1}$ and $\alpha=-\frac{\alpha_2}{\alpha_1}$ with
$\alpha_1\not=0$,
where $\alpha_1$ and $\alpha_2$ are given in {\rm (\ref{lsunified})}.
Suppose $\left(\frac{1}{\lambda-\sigma},x\right)$ is a simple
desired eigenpair of $B$ and
let $(x,X_\perp)$ be unitary. Then
$$
\left[\begin{array}{c}x^H \\ X_\perp^H \end{array}\right]B
\left[\begin{array}{cc}x & X_\perp \end{array}\right]=
\left[\begin{array}{cc}\frac{1}{\lambda-\sigma} & c^H \\
0 & L\end{array}\right],
$$
where $c^H=x^HBX_\perp$ and $L=X_\perp^HBX_\perp$.
Assume that $\sin\angle(\mV,f)\neq0$ and
$\alpha$ is not an eigenvalue of $L$ and define
$$
\sep\left(\alpha ,L\right)=\|(L-\alpha  I)^{-1}\|^{-1}>0.
$$
Then
\begin{equation}\label{epsilonrelation}
\varepsilon \leq C \tilde{\varepsilon},
\end{equation}
where
\begin{equation}\label{defC}
C :=\frac{2\|B\|}{\sep\left(\alpha,L\right)
\sin\angle(\mV,f )}
\end{equation}
with $\alpha=\frac{1}{\nu-\sigma}$ or $\alpha=\frac{1}{\rho-\sigma}$
in SIRA type methods and
$\alpha=y^HBy$ or $\hat{y}^HB\hat{y}$ in JD type methods.
\end{theorem}

If $\alpha$ is a fairly approximation to $\frac{1}{\lambda-\sigma}$
and $\frac{1}{\lambda-\sigma}$ is not close to any eigenvalue of $L$, then
$\sep(\alpha,L)$ is not small and is actually of $O(\|B\|)$ provided that the
eigensystem of $B$ is not ill conditioned. In this case, $C$ is moderate
as $\sin\angle(\mV,f)$ is moderate, as commented previously. For a given
$\tilde\varepsilon$,
we should note that the bigger $\sin\angle(\mV,f)$ is, the bigger $\varepsilon$
is. So it is a lucky event if $\sin\angle(\mV,f)$ is big as it means
that we need to solve inner linear systems with less accuracy $\varepsilon$.

Below we discuss how to determine $\tilde{\varepsilon}$
such that the inexact methods mimic their exact counterparts
very well, that is, each inexact method and its exact counterpart
use comparable or almost the same outer iterations to achieve
the convergence. The following important result is proved in \cite{jiali11},
which forms the basis of our analysis.

\begin{theorem}
Define
$$
\mV_+=\mV\oplus\Span\left\{(I-P_\mV)u\right\},~~~~
\tilde{\mV}_+=\mV\oplus\Span\left\{(I-P_\mV)\tilde{u}\right\},
$$
where $u$ and $\tilde{u}$ are the exact and inexact solutions of the linear system
{\rm (\ref{lsunified})}, respectively, and
\begin{equation}\label{onestep}
\delta=\frac{\sin\angle(\mV_+,x)}{\sin\angle(\mV,x)},\ \ \
\tilde{\delta}=\frac{\sin\angle(\tilde{\mV}_+,x)}{\sin\angle(\mV,x)}.
\end{equation}
Then we have
\begin{equation}\label{deftau}
\frac{|\delta-\tilde{\delta}|}{\delta}\leq2\frac{\tilde{\varepsilon}}{\delta}
:=\tau.
\end{equation}
Furthermore, if $\tau<1$, then
\begin{equation}\label{interval}
\frac{\sin\angle(\tilde{\mV}_+,x)}{\sin\angle(\mV_+,x)}
=\frac{\tilde{\delta}}{\delta}\in[1-\tau,1+\tau].
\end{equation}
\end{theorem}

(\ref{onestep}) measures one step subspace improvements of the exact and inexact
subspace expansions, respectively. (\ref{interval}) indicates that
to make $\sin\angle(\tilde{\mV}_+,x)\approx\sin\angle(\mV_+,x)$,
$\tau$ should be small. Note that the difference of the upper and lower bounds
is $2\tau$. So a very small $\tau$ only improves the bounds
very marginally, and a fairly small $\tau$, e.g.,
$\tau\in [0.001, 0.01]$, is enough since we have
$\frac{\sin\angle(\tilde{\mV}_+,x)}{\sin\angle(\mV_+,x)}\in [0.999,1.001]$
or $[0.99,1.01]$ and the lower and upper bounds differ marginally, which
means that the two subspaces $\tilde{\mV}_+$ and $\mV_+$ are of comparable
or almost the same quality for a fairly small $\tau$ when computing
$(\lambda,x)$. Precisely, with respect to the two subspaces
the approximations of the desired $(\lambda,x)$ obtained by
the same type, i.e., the Rayleigh--Ritz method, its harmonic variant or
the refined harmonic variant, should generally have the comparable or almost
the same accuracy. By (\ref{deftau}), we have
\begin{equation}\label{sufftepsilon}
\tilde{\varepsilon}=\frac{\tau\times\delta}{2}.
\end{equation}

As it has been proved in \cite{jiali11}, the size of $\delta$ crucially
depends on the eigenvalue distribution of $B$. Generally speaking,
the better $\frac{1}{\lambda-\sigma}$ is separated
from the other eigenvalues of $B$, the smaller $\delta$ is.
Conversely, if $\frac{1}{\lambda-\sigma}$ is poorly separated from the others,
$\delta$ may be near to one; see \cite{jiali11} for a descriptive analysis.
$\delta$ is an a-priori quantity and is unknown during computation. But generally
we should not expect that a practical problem is very well conditioned, that is,
$\delta$ is not much less than one;
otherwise, the exact methods will generally find $(\lambda,x)$ by using
only a very few outer iterations. For example, suppose
that the initial $\mV$ is one dimensional and
$\delta=0.1$ at each outer iteration. Then the updated
$\sin\angle(\mV_+,x)$ is no more than $10^{-10}$ after ten outer iterations
and $\mV_+$ is a very accurate subspace for computing $x$.
So in practice, we assume that $\delta$ is not much less than one, say
no smaller than $0.2$. As we have argued, a fairly small $\tau\in [0.001, 0.01]$
should  make the inexact methods behave very like their exact counterparts. As
a result, by (\ref{sufftepsilon}) and the argument that
a fairly small $\tau\in [0.001,0.01]$
is generally enough , to make the inexact methods mimic
their exact counterparts, it is reasonable to take
\begin{equation}\label{pratepsilon}
\tilde\varepsilon\in[10^{-4},10^{-3}],
\end{equation}
which is also suggested in \cite{jiali11} for the standard SIRA and JD methods.
For $\tau\in [0.001,0.01]$, if $\tilde\varepsilon$ defined by
(\ref{pratepsilon}) is unfortunately bigger than that in (\ref{sufftepsilon}),
the the inexact methods may use more outer iterations and may not behave very
like the exact methods. Even so, however, since solving inner linear
systems with $\tilde\varepsilon$ in (\ref{pratepsilon}) is cheaper than with
considerably smaller $\tilde\varepsilon$, the overall efficiency of
the inexact methods may still be improved.

Finally, for the exact methods, let us, qualitatively and concisely,
show which updated subspace $\mV_+$ is better and which method behaves more
favorable. Keep (\ref{defu}) and $y\in {\mV}$ in mind. The (unnormalized)
expansion vector is
$$
(I-P_{\mV})u=(I-P_{\mV})By.
$$
So we actually use add $By$ to $\mV$ to get the expanded subspace $\mV_+$.
It is easily justified that the better
$y$ approximates $x$, the better does $By$
in direction. So for a more accurate approximate eigenvector $y$,
$By$ and thus the expanded subspace $\mV_+$ contain
richer information on the eigenevector $x$.  As we have argued,
since the harmonic Ritz vector is a more reliable and regular approximation to
$x$ while the Ritz vector may not be
for the interior eigenvalue problem, so the exact HSIRA and HJD may expand
the subspaces better and more regularly than the
exact SIRA and JD. Furthermore, since the refined harmonic Ritz vector is
generally more and can be much more accurate than the harmonic Ritz vector,
the exact RHSIRA and RHJD, in turn, generate better subspaces
than the exact HSIRA and HJD at each outer iteration.
As a consequence, HSIRA and HJD are expected to converge more regularly and
use fewer outer iterations than SIRA and JD while RHSIRA and RHJD
may use the fewest outer iterations among all the exact
methods for the interior eigenvalue problem.
For the inexact methods, provided that the selection of $\tilde\varepsilon$ makes
each inexact method mimic its exact counterpart well, the inexact HSIRA, HJD and
RHSIRA, RHJD are advantageous to the inexact SIRA and JD.
Numerical experiments will confirm our expectations.

\section{Restarted algorithms and stopping criteria for
inner solves}\label{sec:restart}

Due to the storage requirement and computational cost,
it is generally necessary to restart the methods to avoid
large steps of outer iterations. We describe the restarted
HSIRA/HJD algorithms as Algorithm~\ref{alg:rehsirajd}
and their refined variants as Algorithm~\ref{alg:rerhsirajd}.

\begin{algorithm}[!htb]
\caption{Restarted HSIRA/HJD algorithm with the fixed target $\sigma$}
\label{alg:rehsirajd}
\begin{algorithmic}
\STATE{Given the target $\sigma$ and a user-prescribed convergence tolerance $tol$,
suppose the columns of $V$ form an orthonormal basis of an initial subspace
$\mV$ and let $\mathrm{M}_{\max}$ be the maximum of outer iterations allowed.}
\REPEAT
\STATE{\begin{enumerate}
\setlength{\itemsep}{0pt}
\setlength{\parsep}{0pt}
\setlength{\parskip}{0pt}
\item Compute (update) $H=V^H(A-\sigma I)^HV$ and $G=V^H(A-\sigma I)^H(A-\sigma I)V$.
\item Let $(\mu,z)$ be an eigenpair of the matrix pencil $(H,G)$,
where $\mu\cong\lambda$.
\item Compute the Rayleigh quotient $\rho=z^HH^Hz+\sigma$ and the harmonic Ritz vector
$y=Vz$.
\item Compute the residual $r=Ay-\rho y$.
\item In HSIRA or HJD, solve the inner linear system
$$
(A-\sigma I)u=r~~~~{\rm or}~~~~(I-yy^H)(A-\sigma I)(I-yy^H)u=-r.
$$
\item Orthonormalize $u$ against $V$ to get the expansion vector $v$.
\item If $\dim(\mathcal{V})< \mathrm{M}_{\max}$,
set $V=\left[\begin{array}{cc}V & v\end{array}\right]$;
otherwise, $V=V_{\mathrm{new}}$.
\end{enumerate}}
\UNTIL{$\|r\|<tol$}
\end{algorithmic}
\end{algorithm}
\begin{algorithm}[!htb]
\caption{Restarted RHSIRA/RHJD algorithm with the fixed target $\sigma$}
\label{alg:rerhsirajd}
\begin{algorithmic}
\STATE{Given the target $\sigma$ and a user-prescribed convergence tolerance $tol$,
suppose the columns of $V$ form an orthonormal basis of an initial subspace
$\mV$ and let $\mathrm{M}_{\max}$ be the maximum of outer iterations allowed.}
\REPEAT
\STATE{
\begin{enumerate}
\setlength{\itemsep}{0pt}
\setlength{\parsep}{0pt}
\setlength{\parskip}{0pt}
\item Compute (update) $H=V^H(A-\sigma I)^HV$ and $G=V^H(A-\sigma I)^H(A-\sigma I)V$.
\item Let $(\mu,z)$ be an eigenpair of the matrix pencil $(H,G)$,
where $\mu\cong\lambda$.
\item Compute the Rayleigh quotient $\rho=z^HH^Hz+\sigma$.
\item Form the cross-product matrix
$S=G+(\overline{\sigma-\rho})H^H+(\sigma-\rho)H+|\sigma-\rho|^2I$ and
compute the eigenvector $\hat{z}$ of $S$ associated with its smallest eigenvalue.
\item Compute the new Rayleigh quotient
$\rho=\hat{z}^HH^H\hat{z}+\sigma$ and the refined harmonic Ritz
vector $y=V\hat{z}$.
\item Compute the residual $r=Ay-\rho y$.
\item In RHSIRA or RHJD, solve the inner linear system
$$
(A-\sigma I)u=r~~~~{\rm or}~~~~(I-yy^H)(A-\sigma I)(I-yy^H)u=-r.
$$
\item Orthonormalize $u$ against $V$ to get the expansion vector $v$.
\item If $\dim(\mathcal{V})< \mathrm{M}_{\max}$, set $V=\left[\begin{array}{cc}V &
v\end{array}\right]$; otherwise, $V=V_{\rm new}$.
\end{enumerate}}
\UNTIL{$\|r\|<tol$}
\end{algorithmic}
\end{algorithm}

We consider some practical issues on two algorithms. For outer iteration steps
$m=1,2,\ldots,\mathrm{M}_{\max}$ during the current cycle,
suppose $(\rho^{(m)},y^{(m)})$ is used to approximate
the desired eigenpair $(\lambda,x)$ of $A$ at the $m$-th outer iteration.
As done commonly in the literature,
we simply take the new starting vector
$$
V_{\mathrm{new}}=y^{(\mathrm{M}_{\max})}
$$
to restart the algorithms. In practice, if $y^{(\mathrm{M}_{\max})}$ is complex
conjugate, we take their real and imaginary parts, normalize and orthonormalize
them to get an orthonormal $V_{\mathrm{new}}$ of column two. We
mention that a thick restarting technique  \cite{stath98} may be used,
in which, besides $y^{(\mathrm{M}_{\max})}$, one retains the approximate
eigenvectors associated with a few other approximate eigenvalues closest to
$\rho^{(\mathrm{M}_{\max})}$ and orthonormalizes them to obtain
$V_{\mathrm{new}}$. We will not report numerical results with thick restarting.

Note that all the methods under consideration do not have
residual monotonically decreasing properties as outer iteration steps increase.
Therefore, it may be possible for both Algorithm \ref{alg:rehsirajd} and
Algorithm \ref{alg:rerhsirajd} to take bad restarting vectors.
This is indeed the case for the standard Rayleigh--Ritz method when computing
interior eigenpairs, as has been widely realized, e.g., \cite{vandervorst2002eigenvalue}.
However, as was stated previously, the harmonic and refined harmonic methods are more
reliable to select correct and good approximations to the desired eigenpairs. So they
are more suitable and reliable for restarting than the standard Ritz vector for
the interior eigenvalue problem. As a consequence, Algorithm~\ref{alg:rehsirajd} and
especially Algorithm~\ref{alg:rerhsirajd} converge more regularly than
the restarted standard SIRA and JD. Our numerical examples will confirm this and
illustrate that the refined harmonic Ritz vectors are better
than the harmonic Ritz vectors for restarting and Algorithm~\ref{alg:rerhsirajd}
converges faster than  Algorithm~\ref{alg:rehsirajd}.

Our ultimate goal is to get a practical estimate for $\varepsilon$ for a
given $\tilde{\varepsilon}$. The a-priori bound (\ref{epsilonrelation}) forms
the basis for this. Let $\nu_1,\nu_2,\ldots,\nu_m$ are the harmonic Ritz values,
i.e., the eigenvalues of the pencil $(H_m,G_m)$ at the $m$-th outer iteration during the
current cycle (here we add the subscript $m$ to $H$ and $G$ in the algorithms) and
assume that $\nu_1$ is used to approximate the desired eigenvalue $\lambda$. We simply
estimate $\|B\|\approx \frac{1}{\mid\rho-\sigma\mid}$, where $\rho$ is the
Rayleigh quotient of $A$ with respect to the harmonic Ritz vector in
HSIRA and HJD and the refined harmonic Ritz vector in RHSIRA and RHJD. For
$\sep (\alpha,L)$, we replace $\alpha$, an approximation to
$\frac{1}{\lambda-\sigma}$, by $\frac{1}{\rho-\sigma}$ accordingly and
then estimate
$$
\sep (\alpha,L)\approx \min_{i=2,3,\ldots,m}
\mid\frac{1}{\rho-\sigma}-\frac{1}{\nu_i-\sigma}\mid.
$$
Finally, replace $\sin\angle(\mV,f)$ by its maximum one. Combining all these together,
we get the following estimate $C'$ of $C$ in (\ref{epsilonrelation}):
\begin{equation}\label{Cstar2}
C'=
\left\{\begin{array}{rl}
2\max\limits_{i=2,3,\ldots,m}\left|\frac{\nu_i-\sigma}{\nu_i-\rho}\right|,&m>1\\
1,&m=1
\end{array}\right.,
\end{equation}
This is analogous to what is done in \cite{jiali11}, and the difference
is that we here use the harmonic Ritz values
$\nu_2,\ldots,\nu_m$ to replace the standard Ritz values used in \cite{jiali11}
as approximations to some eigenvalues of $A$ other than $\lambda$.
Denote by $\varepsilon_S$ and $\varepsilon_J$ practical $\varepsilon$'s used
in the HSIRA, RHSIRA and HJD, RHJD algorithms, respectively. Recall
that bound (\ref{epsilonrelation}) is compact. Then we take
\begin{equation}\label{epsilon'}
\varepsilon_S=\varepsilon_J=\varepsilon=C'\tilde{\varepsilon}.
\end{equation}
For a fairly small $\tilde\varepsilon$, we may have $\varepsilon\geq1$ in case $C'$
is big, which will make $\tilde{u}$ no sense as an approximation to $u$.
In order to make $\tilde{u}$ have some accuracy,
we propose using
\begin{equation}\label{epsmin}
\varepsilon=\min\left\{C'\tilde{\varepsilon},0.1\right\},
\end{equation}
which shows that $\varepsilon$ is comparable to $\tilde{\varepsilon}$ in size
whenever $\tilde\varepsilon$ is small and $C'$ is moderate.

\section{Numerical experiments}\label{sec:experiments}

All the numerical experiments were performed on an Intel
(R) Core (TM)2 Quad CPU Q9400 $2.66$GHz with the main memory 2 GB using
Matlab 7.8.0 with the machine precision $\epsilon_{\rm mach}=2.22\times
10^{-16}$ under the Linux operating system.
All the test examples are difficult in the sense that the desired eigenvalue
is clustered with some other eigenvalues of $A$
for the given target $\sigma$. We aim to show four points.
First, regarding restarts of outer iterations,
for fairly small $\tilde{\varepsilon}=10^{-3}$
and $10^{-4}$, the restarted inexact HSIRA/HJD and RHSIRA/RHJD algorithms
behave (very) like the exact counterparts.
Second, regarding outer iterations,
Algorithms~\ref{alg:rehsirajd}--\ref{alg:rerhsirajd} are much more efficient
than the restarted standard SIRA and JD for the same $\tilde\varepsilon$'s, and
Algorithm~\ref{alg:rerhsirajd} is the best.
Third, regarding total inner iterations and overall efficiency,
Algorithm~\ref{alg:rerhsirajd} is considerably more efficient than
Algorithm~\ref{alg:rehsirajd}, and the restarted
inexact standard SIRA and JD perform very poorly and often fail to converge.
Fourth, each SIRA type algorithm is equally as efficient as the corresponding
JD type algorithm for the same $\tilde\varepsilon$.

At the $m$-th outer iteration step of the inexact HSIRA or HJD method, we have
$H_m=V_m^H(A-\sigma I)^HV_m$ and $G_m=V_m^H(A-\sigma I)^H(A-\sigma I)V_m$.
Let $(\nu_i^{(m)},V_mz_i^{(m)}),\ i=1,2,\ldots,m$ be the harmonic Ritz pairs,
labeled as $|\nu_1^{(m)}-\sigma|<|\nu_2^{(m)}-\sigma|\leq\cdots\leq
|\nu_m^{(m)}-\sigma|$. Keep (\ref{rho}) in mind. We use
\begin{equation}\label{approxpair}
(\rho^{(m)},y^{(m)})=\left((z_1^{(m)})^HH_m^Hz_1^{(m)}+\sigma,V_mz_1^{(m)}\right)
\end{equation}
to approximate the desired eigenpair $(\lambda,x)$ of $A$.
If the RSIRA and RHJD methods are used, we form the cross-product matrix
$S_m$ by (\ref{comS}) and compute its eigenvector $\hat{z}^{(m)}$ associated with
the smallest eigenvalue. We then compute the refined harmonic Ritz vector
$y^{(m)}=V_m\hat z^{(m)}$ and the Rayleigh quotient $\rho^{(m)}$ defined
by (\ref{uprho}). Let $r_m=Ay^{(m)}-\rho^{(m)}y^{(m)}$. We stop the algorithms if
\begin{eqnarray}
\|r_m\|< tol=\max\left\{\|A\|_1,1\right\}\times10^{-12}.
\end{eqnarray}

In the algorithms, we stop inner iterations at
outer iteration $m$ when
\begin{equation}\label{stopinner}
\frac{\|r_m-(A-\sigma I)\tilde{u}\|}{\|r_m\|},~~~
\frac{\|-r_m-(I-y^{(m)}(y^{(m)})^H)(A-\sigma I)(I-y^{(m)}(y^{(m)})^H)\tilde{u}\|}
{\|r_m\|}\leq\varepsilon,
\end{equation}
where $\varepsilon$ is defined by (\ref{epsmin}) for a given $\tilde\varepsilon$.
We will denote by SIRA($\tilde\varepsilon$), JD($\tilde\varepsilon$),
HSIRA($\tilde\varepsilon$), HJD($\tilde\varepsilon$) and RHSIRA($\tilde\varepsilon$),
RHJD($\tilde\varepsilon$) the inexact algorithms with the given $\tilde\varepsilon$.

In the \textquotedblleft exact\textquotedblright\ SIRA and JD type
algorithms, we stop inner iterations when (\ref{stopinner}) is satisfied
with $\varepsilon=10^{-14}$.

All the test problems are from Matrix Market \cite{matrixmarket}.
For each inner linear system, we always took the zero vector as
an initial approximate solution and solved it
by the right-preconditioned restarted GMRES(30) algorithm.
Each algorithm starts with the normalized vector
$\frac{1}{\sqrt{n}}(1,1,\ldots,1)^H$.
At the $m$-th outer iteration of current restart, $m=1,2,\ldots,{\rm M}_{\max}$,
for the correction equations in the JD type algorithms we used the preconditioner
\begin{equation}\label{precondcor}
M_m=\left(I-y^{(m)}(y^{(m)})^H\right)M\left(I-y^{(m)}(y^{(m)})^H\right)
\end{equation}
suggested in \cite{vandervorst2002eigenvalue},
where $M\approx A-\sigma I$ is the incomplete LU preconditioner
for the SIRA type algorithms.
For all the algorithms, the maximum steps of outer iterations
are $30$ per restart. An algorithm signals the failure if it did not converge
within ${\bf Max}$ restarts.

In the experiments, each figure consists of three subfigures, in which
the top subfigure denotes outer residual norms versus outer iterations
of the first cycle before restart, the bottom-left subfigure denotes outer
residual norms versus restarts and the bottom-right subfigure denotes
the numbers of inner iterations versus restarts.

The top subfigures exhibit the qualities of
the standard, harmonic and refined harmonic Ritz vectors for restarting.
The bottom two subfigures depict the convergence processes of
the exact SIRA type algorithms and their inexact counterparts
with $\tilde{\varepsilon}=10^{-3}$, showing the convergence behavior of
the algorithms and the local efficiency per restart, respectively.

In all the table, denote by $I_{restart}$ the number of restarts to achieve the
convergence, by $I_{inner}$ the total number of inner iterations and by
$P_{0.1}$ the percentage of the times $I_{0.1}$ that $\varepsilon=0.1$
is used in the total number of outer iterations $I_{outer}$, i.e.,
$$
P_{0.1}=\frac{I_{0.1}}{I_{outer}}.
$$
Note that $I_{inner}$ equals the total products of $A$ and vectors
in the restarted GMRES algorithm. It is a reasonable measure of
the overall efficiency of all the algorithms used in the experiments.
\smallskip

\begin{figure}[!htb]
\centering
\includegraphics[width=13cm,height=5.5cm]{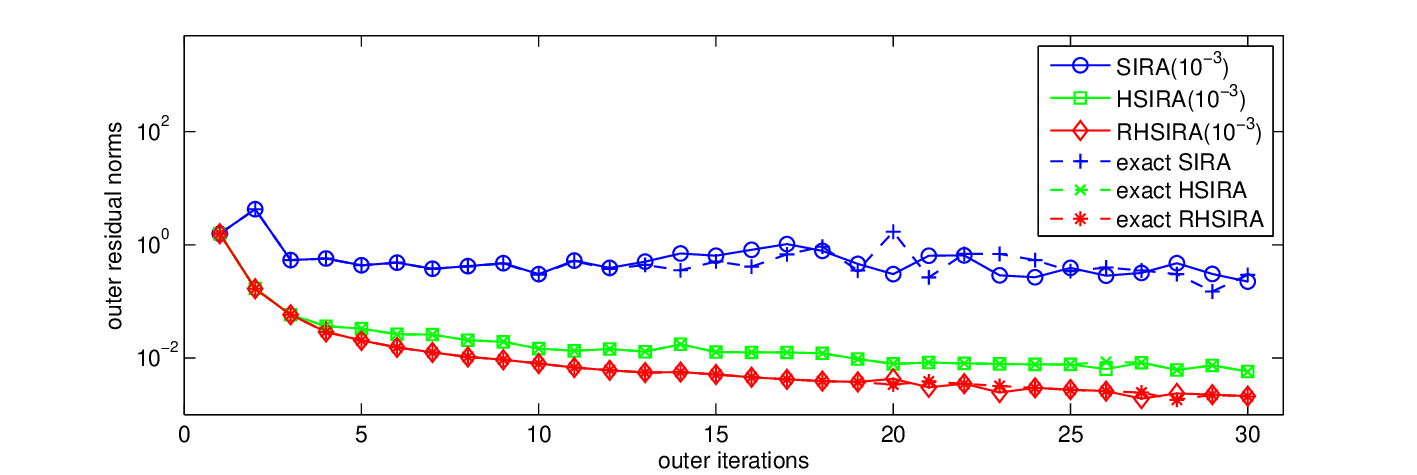}
\begin{minipage}{7.5cm}
\includegraphics[width=8.3cm,height=5.5cm]{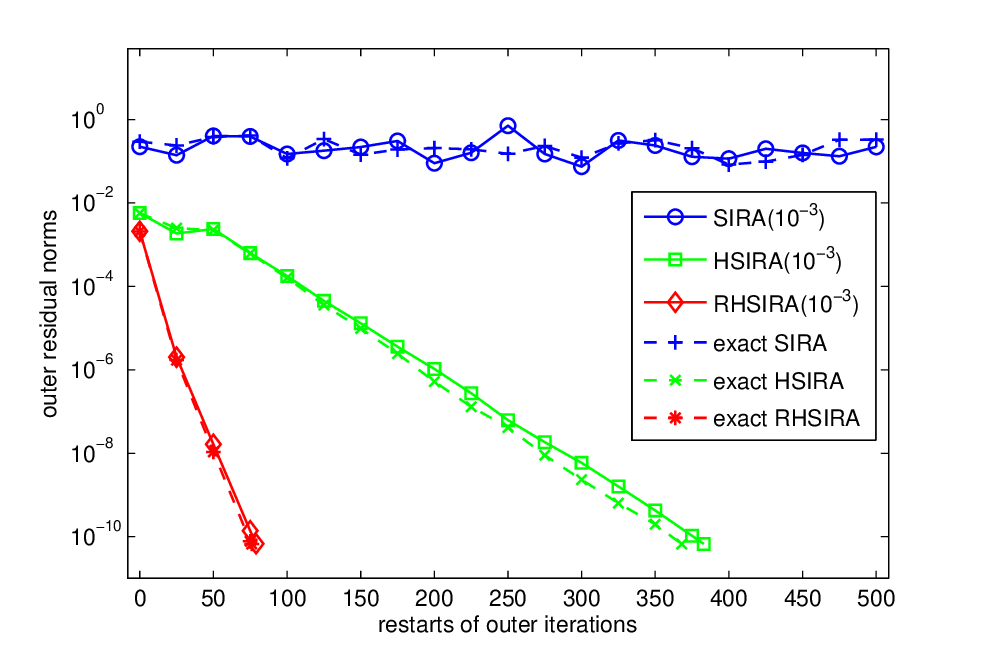}
\end{minipage}
\begin{minipage}{7.5cm}
\includegraphics[width=8.3cm,height=5.5cm]{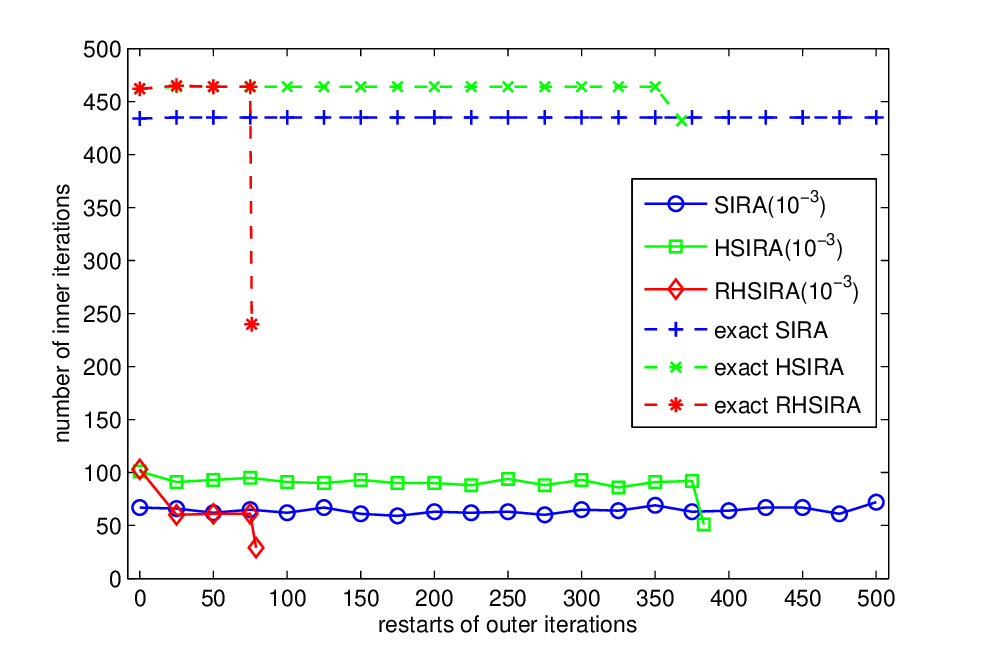}
\end{minipage}
\caption{\emph{Example 1. M80PI\_n with $\sigma=0.05+0.5\imag$.}}
\label{fig_m80pi_n}
\end{figure}

\noindent\textbf{Example 1.} This unsymmetric eigenproblem M80PI\_n of $n=4182$
arises from real power system models \cite{matrixmarket}.
We test this example with $\sigma=0.05+0.5\imag$ and ${\bf Max}=500$.
The computed eigenvalue is
$\lambda\approx -6.9391\times10^{-5}+5.0062\times10^{-1}\imag$.
The preconditioner $M$ is the incomplete LU factorization of
$A-\sigma I$ with drop tolerance $10^{-1}$. Figure~\ref{fig_m80pi_n} and
Table~\ref{tab_m80pi_n} display the results.
\smallskip

\begin{table}[!htb]
\centering
\small
\begin{tabular}{c||c|r|r|r||c|r|r|r}
\hline
Accuracy &  Algorithm  &  $I_{restart}$  & $I_{inner}$ & $P_{0.1}$ &
Algorithm  &  $I_{restart}$  & $I_{inner}$ & $P_{0.1}$\\
\hline\hline
& SIRA & ${\bf Max}$ & $32142$ & $0\%$
& JD & ${\bf Max}$   & $33823$ & $0\%$   \\
$\tilde{\varepsilon}=10^{-3}$
& HSIRA & $383$   & $35032$ & $0\%$
& HJD & $371$   & $33988$ & $0\%$ \\
& RHSIRA & $79$  & $4981$ & $89\%$
& RHJD & $79$  & $4958$ & $89\%$\\
\hline
& SIRA & ${\bf Max}$   & $47166$ & $0\%$
& JD & ${\bf Max}$   & $47785$ & $0\%$    \\
$\tilde{\varepsilon}=10^{-4}$
& HSIRA & $376$   & $46470$ & $0\%$
& HJD & $374$   & $46231$ & $0\%$    \\
& RHSIRA & $74$  & $7032$ & $0\%$
& RHJD & $75$  & $7148$  & $0\%$  \\
\hline
& SIRA & ${\bf Max}$  & $217941$ & $-$
& JD & ${\bf Max}$  & $217936$ & $-$   \\
exact
& HSIRA & $368$   & $171183$  & $-$
& HJD & $368$   & $171182$ & $-$      \\
& RHSIRA & $76$  & $35515$ & $-$
& RHJD & $76$  & $35513$ & $-$ \\
\hline
\end{tabular}
\caption{\emph{Example 1. M80PI\_n with $\sigma=0.05+0.5\imag$.}}
\label{tab_m80pi_n}
\end{table}

We see from the top subfigure of Figure~\ref{fig_m80pi_n}
that for $\tilde\varepsilon=10^{-3}$ the inexact SIRA, HSIRA and RHSIRA
behaved very like their corresponding exact counterparts in the first cycle
but the exact and inexact RHSIRA are more effective than the exact
and inexact HSIRA. It clearly shows that the exact SIRA and SIRA($10^{-3}$)
are the poorest and considerably poorer than the corresponding HSIRA and RHSIRA.
There are two reasons for this. The first is that
harmonic and refined harmonic Ritz vectors are more reliable
than standard Ritz vectors for computing interior eigenvectors.
The second is that the harmonic and refined  harmonic Ritz vectors
favor subspace expansions. We, therefore,
expect that the restarted HSIRA and RHSIRA can be much more efficient and converge
much faster than the restarted SIRA when computing interior eigenpairs
since we use more reliable and possibly accurate harmonic and refined
harmonic Ritz vectors as restarting vectors. Furthermore,
as far as restarts are concerned, the restarted RHSIRA outperforms the restarted
HSIRA very substantially. We observed very similar convergence behavior for the JD
type algorithms and thus have the same expectations on the restarted JD type
algorithms. These expectations are indeed confirmed by numerical experiments, as
shown by $I_{restart}$'s in the table and figure.

We explain the table and figure in more details. From the bottom-left subfigure,
regarding outer iterations, we see that both the exact and inexact
restarted SIRA did not converge within $500$ restarts while HSIRA and RHSIRA
worked well and the inexact methods behaved very like their
corresponding exact ones. The restarted RHSIRA and RHJD were the
fastest and five times as fast as the restarted HSIRA and HJD, respectively.
As is expected, the table confirms that, for the same $\tilde\varepsilon$,
each SIRA type algorithm and the corresponding JD type
algorithm behaved very similar and were almost
indistinguishable.

Regarding the overall efficiency, the exact HSIRA, HJD, RHSIRA and RHJD each
used about $460$ inner iterations per restart. In contrast, HSIRA($10^{-3})$ and
HSIRA($10^{-4}$) used almost constant inner iterations each restart, which were
about $100$ and 110 per restart, respectively, and
RHSIRA($10^{-3}$) and RHSIRA($10^{-4}$) used roughly $60$ and 100 inner iterations
per restart, respectively. The same observations are true for the JD type algorithms.
These experiments demonstrate that modest $\tilde\varepsilon=10^{-3}, \ 10^{-4}$
are enough to make the inexact restarted algorithms mimic their exact counterparts.
We find that the SIRA($10^{-4}$) and JD($10^{-4}$) type algorithms
used almost the same outer iterations as the SIRA($10^{-3}$) and JD($10^{-3}$)
type ones, but the latter consumed considerably fewer total inner
iterations and improved the overall efficiency substantially.
So, smaller $\tilde\varepsilon$'s are not necessary for this example as they cannot
reduce outer iterations further and may cost much more inner iterations.

In addition, we see from Table~\ref{tab_m80pi_n} that $89$ percent
of inner linear systems in RHSIRA($10^{-3}$) were actually solved with
the accuracy requirement $\varepsilon=0.1$. This means that though most of $C'$'s in
(\ref{Cstar2}) are big, it suffices to solve all the inner linear
systems with low accuracy $0.1$.
\smallskip

\begin{figure}[!htb]
\centering
\includegraphics[width=13cm,height=5.5cm]{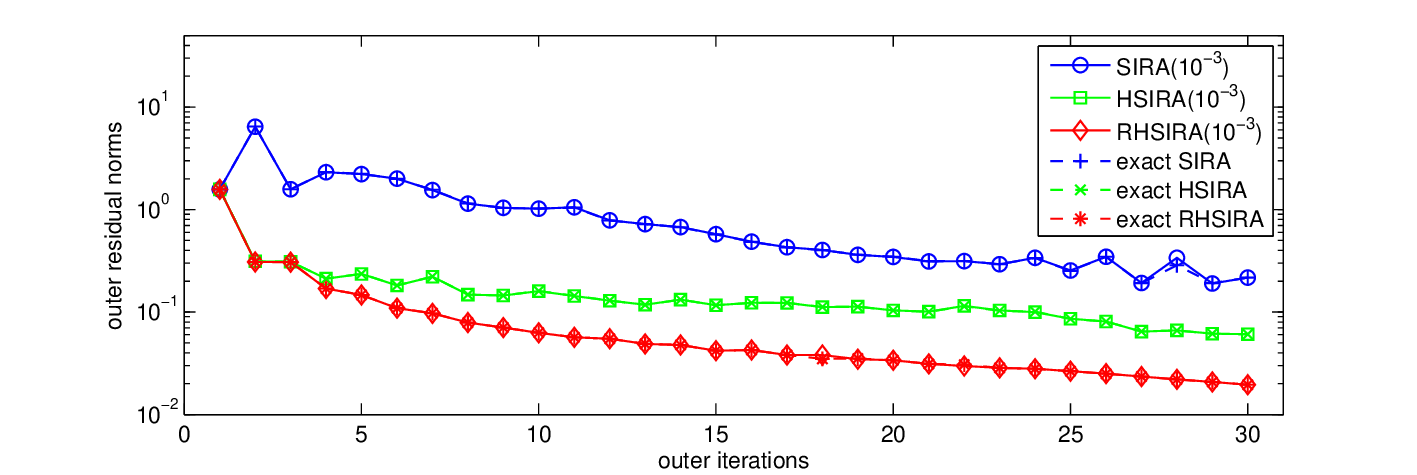}
\begin{minipage}{7.5cm}
\includegraphics[width=8.3cm,height=5.5cm]{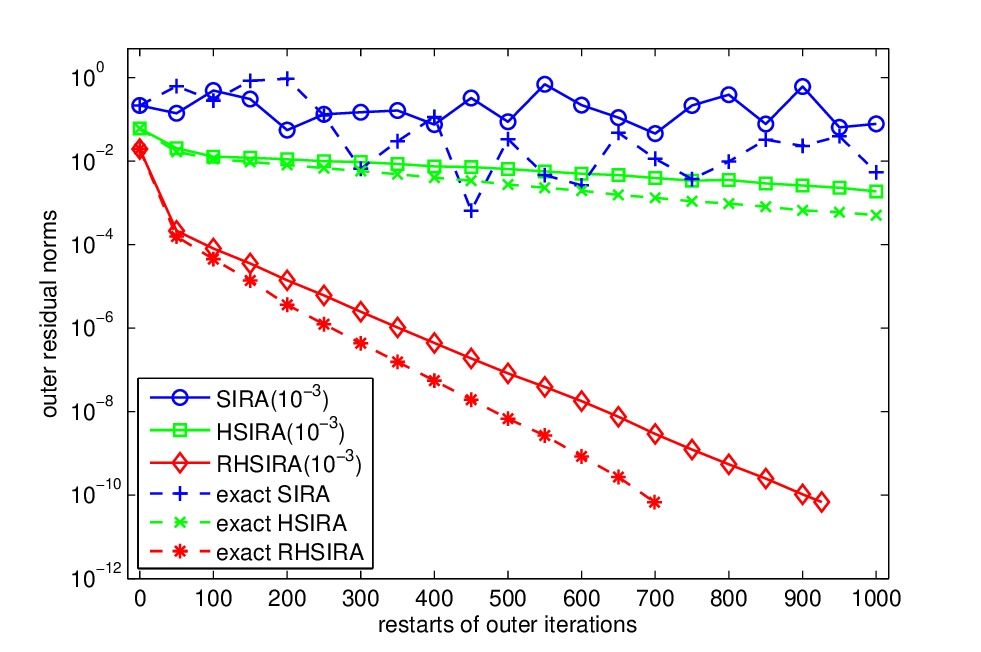}
\end{minipage}
\begin{minipage}{7.5cm}
\includegraphics[width=8.3cm,height=5.5cm]{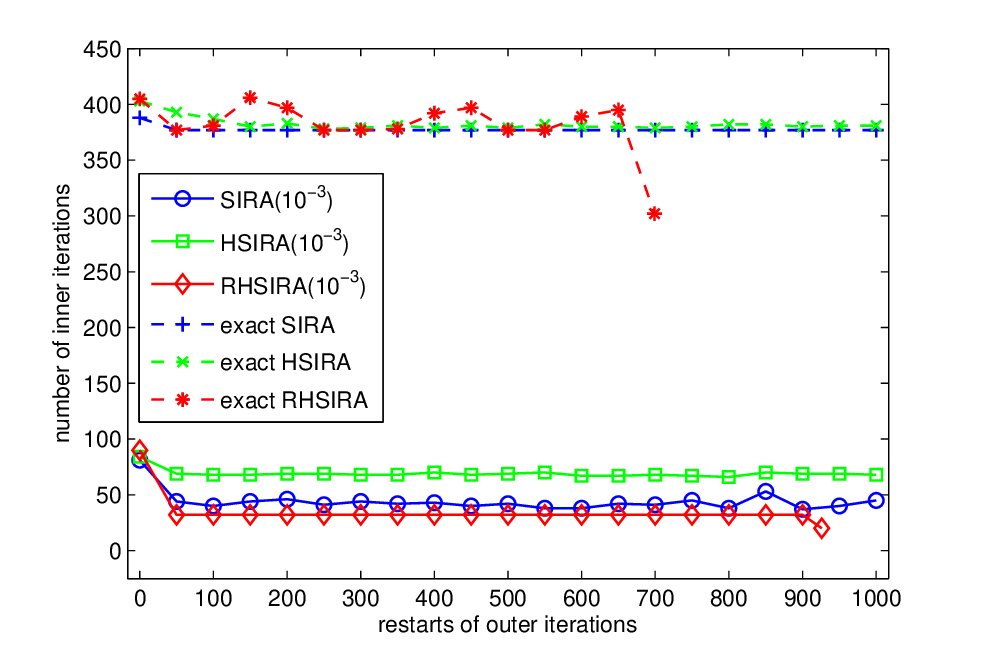}
\end{minipage}
\caption{\emph{Example 2. S80PI\_n with $\sigma=0.4+1.3\imag$.}}
\label{fig_s80pi_n}
\end{figure}

\noindent\textbf{Example 2.} This unsymmetric eigenproblem M80PI\_n of $n=4182$
arises from real power system models \cite{matrixmarket}. We test test this
example with $\sigma=0.4+1.3\imag$ and ${\bf Max}=1000$. The computed eigenvalue is
$\lambda\approx -2.8266\times10^{-4}+1.3068\imag$.
The preconditioner $M$ is the incomplete LU factorization of
$A-\sigma I$ with drop tolerance $10^{-1}$. Figure~\ref{fig_s80pi_n} and
Table~\ref{tab_s80pi_n} report the results.

\begin{table}[!htb]
\centering
\small
\begin{tabular}{c||c|r|r|r||c|r|r|r}
\hline
Accuracy &  Algorithm  &  $I_{restart}$  & $I_{inner}$ & $P_{0.1}$ &
Algorithm  &  $I_{restart}$  & $I_{inner}$ & $P_{0.1}$\\
\hline\hline
& SIRA & ${\bf Max}$ & $43020$ & $6\%$
& JD & ${\bf Max}$   & $43279$ & $5\%$   \\
$\tilde{\varepsilon}=10^{-3}$
& HSIRA & ${\bf Max}$   & $67645$ & $0\%$
& HJD & ${\bf Max}$   & $67686$ & $0\%$ \\
& RHSIRA & $926$  & $29774$ & $96\%$
& RHJD & $946$  & $30392$ & $96\%$\\
\hline
& SIRA & ${\bf Max}$   & $82712$ & $0\%$
& JD & ${\bf Max}$   & $82335$ & $0\%$    \\
$\tilde{\varepsilon}=10^{-4}$
& HSIRA & ${\bf Max}$   & $104141$ & $0\%$
& HJD & ${\bf Max}$   & $104102$ & $0\%$    \\
& RHSIRA & $945$  & $31323$ & $94\%$
& RHJD & $937$  & $31068$  & $93\%$  \\
\hline
& SIRA & ${\bf Max}$  & $377424$ & $-$
& JD & ${\bf Max}$  & $377426$ & $-$   \\
exact
& HSIRA & ${\bf Max}$   & $382069$  & $-$
& HJD & ${\bf Max}$   & $381359$ & $-$      \\
& RHSIRA & $699$  & $271142$ & $-$
& RHJD & $684$  & $263865$ & $-$ \\
\hline
\end{tabular}
\caption{\emph{Example 2. S80PI\_n with $\sigma=0.4+1.3\imag$.}}
\label{tab_s80pi_n}
\end{table}

The top subfigure of Figure~\ref{fig_s80pi_n} is similar to that of
Figure~\ref{fig_m80pi_n}, showing that harmonic and especially
refined harmonic Ritz vectors are more suitable for expanding subspaces
and restarting for the  interior eigenvalue problem. The difference is that this
problem is more difficult than Example 1 since, for the
30-dimensional subspace in the first cycle, the residual norms
decrease more slowly and approximate eigenpairs are less accurate than
those for Example 1. So it is expected that the restarted algorithms
converge more slowly and use more outer iterations $I_{restart}$'s
than for Example 1. We observe from the figure
that the convergence curves of the three exact algorithms SIRA, HSIRA and RHSIRA
essentially coincide with those of their inexact variants with $\varepsilon=10^{-3}$.
So smaller $\tilde\varepsilon$ doe not help, rather it makes the algorithms
may waste much more inner iterations. We observed similar convergence phenomena for
the JD type algorithms.

Table~\ref{tab_s80pi_n} and Figure~\ref{fig_s80pi_n} shows that the restarted exact
and inexact SIRA, JD, HSIRA and HJD all failed to converge within 1000 restarts.
Furthermore,  it is deduced from Figure~\ref{fig_s80pi_n} that the restarted exact
and inexact SIRA and JD appear impossible to converge at all since their convergence
curves were irregular, oscillated and had no decreasing tendency. The restarted exact
and inexact SIRA and JD behaved regular but converged too slowly. In contrast,
the restarted exact and inexact RHSIRA and RHJD with $\tilde\varepsilon=10^{-3}$ and
$10^{-4}$ converged. Furthermore, the restarted RHSIRA and RHJD with $\tilde\varepsilon
=10^{-3}$ and $\tilde\varepsilon=10^{-4}$ behaved very similar and used almost
the same outer iterations. Both of them mimic the restarted exact RHSIRA and RHJD fairly
well. So, smaller $\varepsilon$ cannot reduce outer iterations
substantially.

Regarding the overall efficiency, the inexact restarted algorithms improved the
situation tremendously. As in Example 1, the exact RHSIRA and RHJD still
needed about 370 inner iterations per restart, much more than 40 inner
iterations that were used by the inexact RHSIRA and RHJD for
$\tilde\varepsilon=10^{-3}$ and $10^{-4}$ per
restart. As a whole, $I_{inner}$'s illustrate that the inexact algorithms with
these two $\tilde\varepsilon$'s were about eight times faster than the exact
algorithms, a striking improvement. Besides, note that the number of inner
linear systems that were solved with lower accuracy $\varepsilon=0.1$ in
RHSIRA($10^{-3}$) were $96\%$, while those in HSIRA($10^{-3}$) and
SIRA($10^{-3}$) were $0\%$ and $6\%$, respectively.
This is why RHSIRA($10^{-3}$) used fewer inner iterations than SIRA($10^{-3}$)
and HSIRA($10^{-3}$) per restart, as shown in the bottom-right subfigure.

Although the restarted exact SIRA and HSIRA both failed for this example,
the reasons may be completely different.
In the top subfigure, we find that the convergence curves of SIRA bulged
at the last few steps ($25\sim30$). This is harmful to restarting as
it is very possible to get an unsatisfying restarting
vector once the method bulged at the very last step.
In the bottom-left subfigure, it is seen that the convergence curves of
restarted exact and inexact SIRA were irregular while the convergence curves of HSIRA
decreased smoothly though the algorithm converged quite slow.
So we can infer that HSIRA did not take bad restarting vectors. Remarkably,
the figure and table tell us that RHSIRA converged much faster.
This is due to the fact that the refined harmonic Ritz vector $\hat y$ can be
much more accurate than the corresponding harmonic Ritz vector $y$, so that
restarting vectors were better and the subspaces generated were
more accurate as restarts proceeded.
\smallskip

\begin{figure}[!htb]
\centering
\includegraphics[width=13cm,height=5.5cm]{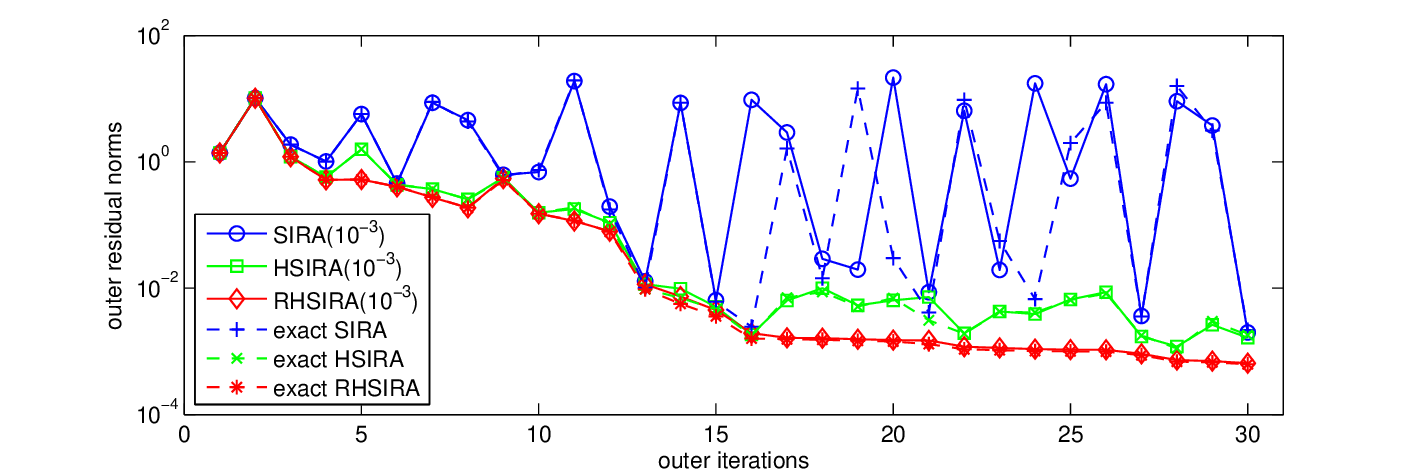}
\begin{minipage}{7.5cm}
\includegraphics[width=8.3cm,height=5.5cm]{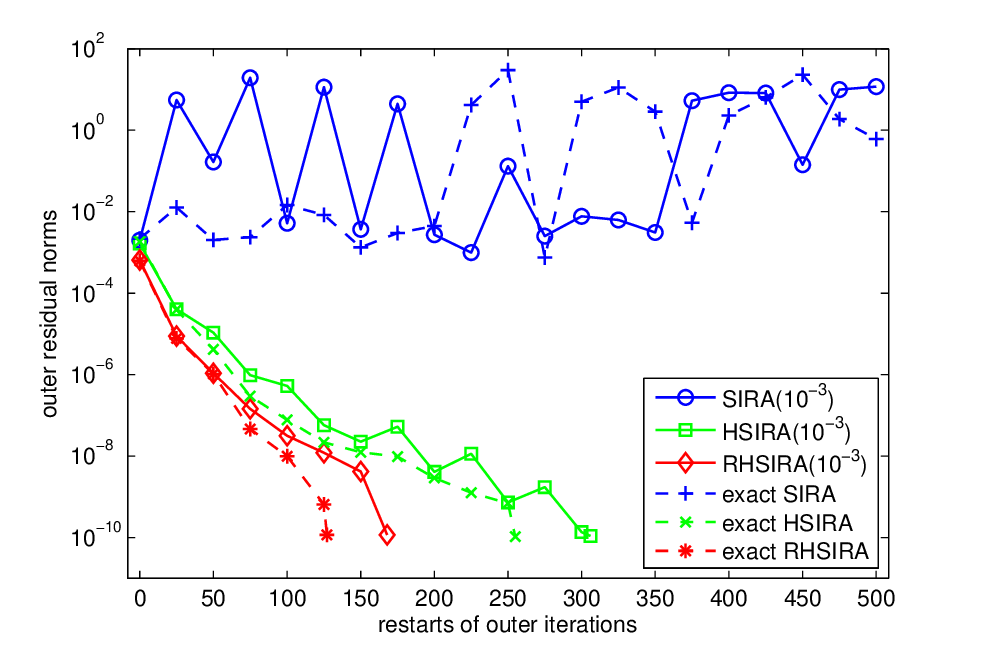}
\end{minipage}
\begin{minipage}{7.5cm}
\includegraphics[width=8.3cm,height=5.5cm]{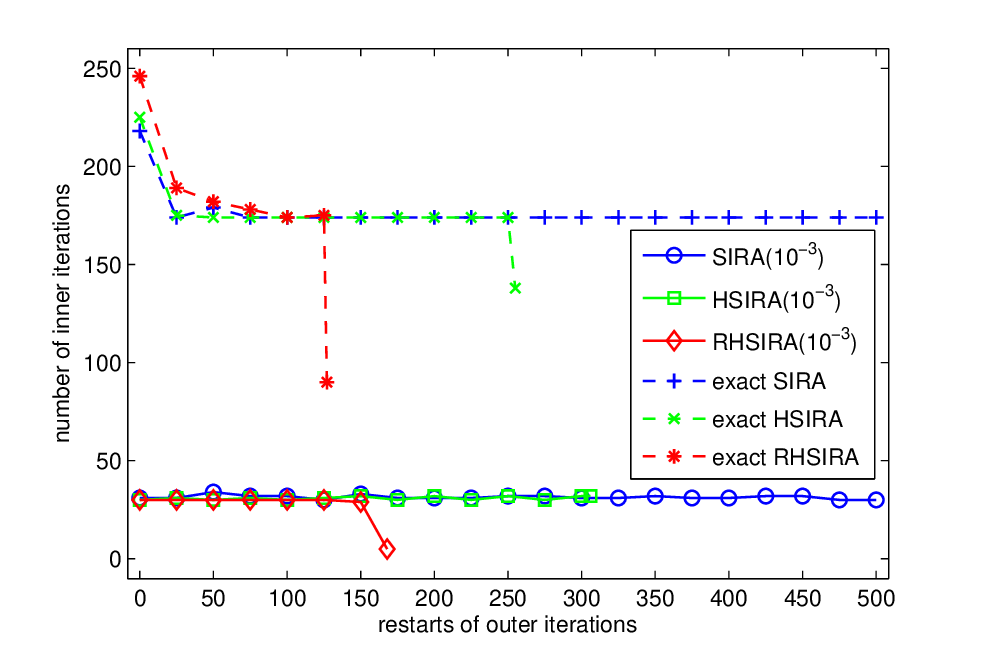}
\end{minipage}
\caption{\emph{Example 3. dw4096 with $\sigma=-24$.}}
\label{fig_dw4096}
\end{figure}

\noindent\textbf{Example 3.} This unsymmetric eigenproblem dw4096 of $n=8192$
arises from dielectric channel waveguide problems \cite{matrixmarket}.
We test this example with $\sigma=-24$ and ${\bf Max}=500$.
The computed eigenvalue is $\lambda\approx-30.217$.
The preconditioner $M$ is the incomplete LU factorization of
$A-\sigma I$ with drop tolerance $10^{-3}$.
Figure~\ref{fig_dw4096} and Table~\ref{tab_dw4096} display the results.
\smallskip

\begin{table}[!htb]
\centering
\small
\begin{tabular}{c||c|r|r|r||c|r|r|r}
\hline
Accuracy &  Algorithm  &  $I_{restart}$  & $I_{inner}$ & $P_{0.1}$ &
Algorithm  &  $I_{restart}$  & $I_{inner}$ & $P_{0.1}$\\
\hline\hline
& SIRA & ${\bf Max}$   & $15722$ & $45\%$
& JD & ${\bf Max}$   & $17118$ & $45\%$   \\
$\tilde{\varepsilon}=10^{-3}$
& HSIRA & $306$   & $9473$ & $85\%$
& HJD & $300$   & $9246$ & $85\%$ \\
& RHSIRA & $168$  & $5019$ & $96\%$
& RHJD & $151$  & $4540$ & $96\%$\\
\hline
& SIRA & ${\bf Max}$   & $20586$ & $14\%$
& JD & ${\bf Max}$   & $20794$ & $14\%$    \\
$\tilde{\varepsilon}=10^{-4}$
& HSIRA & $289$   & $9251$ & $49\%$
& HJD & $282$   & $9091$ & $49\%$    \\
& RHSIRA & $156$  & $4722$ & $94\%$
& RHJD & $197$  & $5952$  & $95\%$  \\
\hline
& SIRA & ${\bf Max}$  & $87410$ & $-$
& JD & ${\bf Max}$  & $89435$ & $-$   \\
exact
& HSIRA & $255$   & $44906$  & $-$
& HJD & $259$   & $45480$ & $-$      \\
& RHSIRA & $127$  & $23351$ & $-$
& RHJD & $139$  & $25419$ & $-$ \\
\hline
\end{tabular}
\caption{\emph{Example 3. dw4096 with $\sigma=-24$.}}
\label{tab_dw4096}
\end{table}

Figure~\ref{fig_dw4096} indicates that the exact and inexact SIRA oscillated sharply
in the first cycle, the exact and inexact HSIRA improved the situation very significantly
but still did not behave quite regularly, while the exact and inexact RHSIRA converged
very smoothly after first a few outer iterations and considerably faster than the
HSIRA. Poor Ritz vectors further led to poor subspace expansion vectors, generating
a sequence of poor subspaces. We observed very similar phenomena in
the JD type algorithms. This means that the standard Ritz vectors were not suitable for
restarting, and the harmonic Ritz vectors were much better but inferior to
the refined harmonic Ritz vectors for restarting. So it is expected that the restarted
exact and inexact SIRA and JD algorithms may not work well, but the restarted
HSIRA and HJD may work much better than the former ones and the restarted RHSIRA and
RHJD outperform the HSIRA and HJD very considerably.

The above expectations are confirmed by  and Table~\ref{tab_dw4096} and
the bottom-left subfigure of Figure~\ref{fig_dw4096}. It is seen that
the restarted exact and inexact SIRA and JD algorithms failed to converge
within 500 restarts while HSIRA and HJD solved the problem successfully
and the refined RHSIRA and RHJD were twice as efficient as the harmonic
algorithms, as indicated by $I_{restart}$'s. Furthermore, we see
the inexact algorithms with $\tilde\varepsilon=10^{-3},\ 10^{-4}$ exhibited
very similar convergence behavior and mimic the exact algorithms very well.
This confirms our theory that modest $\tilde\varepsilon$ is generally
enough to make the inexact SIRA and JD type algorithms behave very like
the corresponding exact counterparts and smaller $\tilde\varepsilon$ is not
necessary. We also find that each exact SIRA type algorithm is as efficient
as the corresponding JD type algorithm for the same $\tilde\varepsilon$.

Regarding the overall performance, the comments on Examples 1--2 apply here
analogously. For $\tilde\varepsilon=10^{-3}$ and $10^{-4}$,
the inexact HSIRA, RHSIRA, HJD and RHJD were four times faster than the
corresponding exact algorithms. Furthermore, the restarted exact and inexact
RHSIRA and RHJD were twice as fast as the corresponding HSIRA and HJD algorithms
for the same $\tilde\varepsilon$.
\smallskip

\begin{figure}[!htb]
\centering
\includegraphics[width=13cm,height=5.5cm]{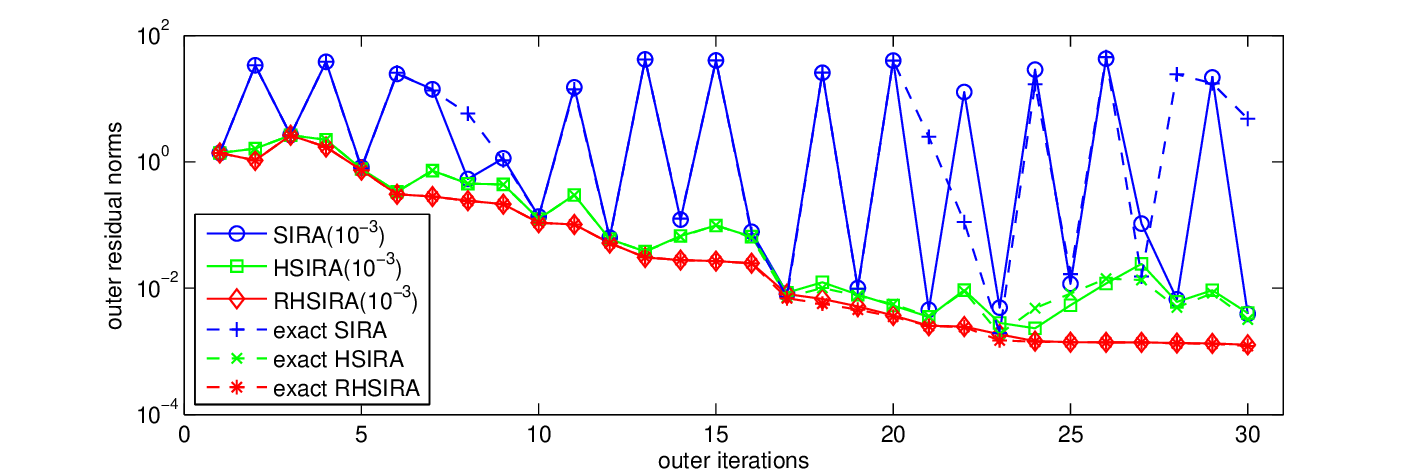}
\begin{minipage}{7.5cm}
\includegraphics[width=8.3cm,height=5.5cm]{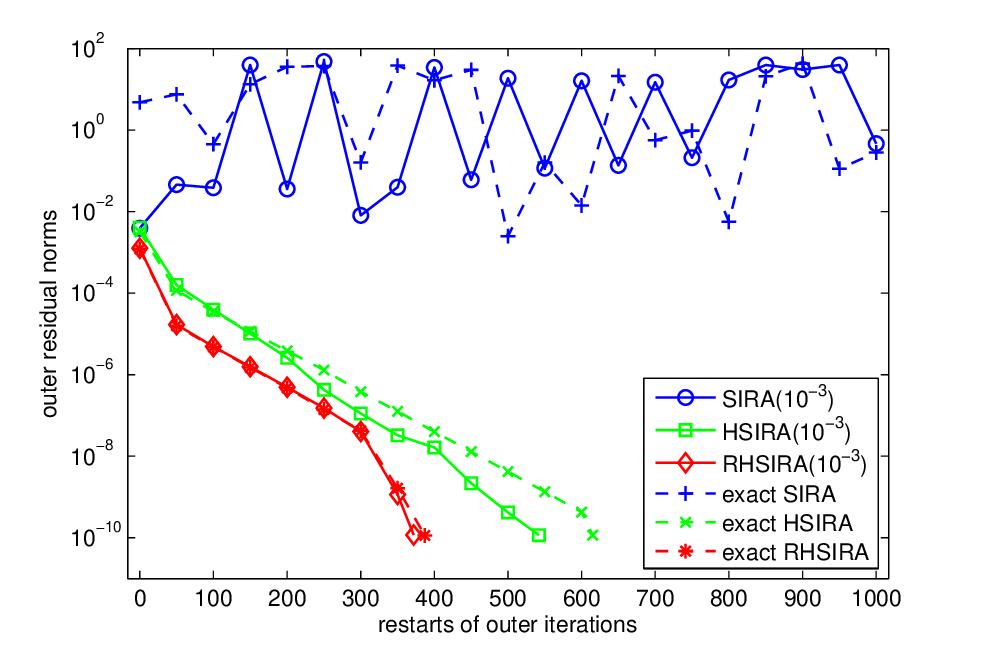}
\end{minipage}
\begin{minipage}{7.5cm}
\includegraphics[width=8.3cm,height=5.5cm]{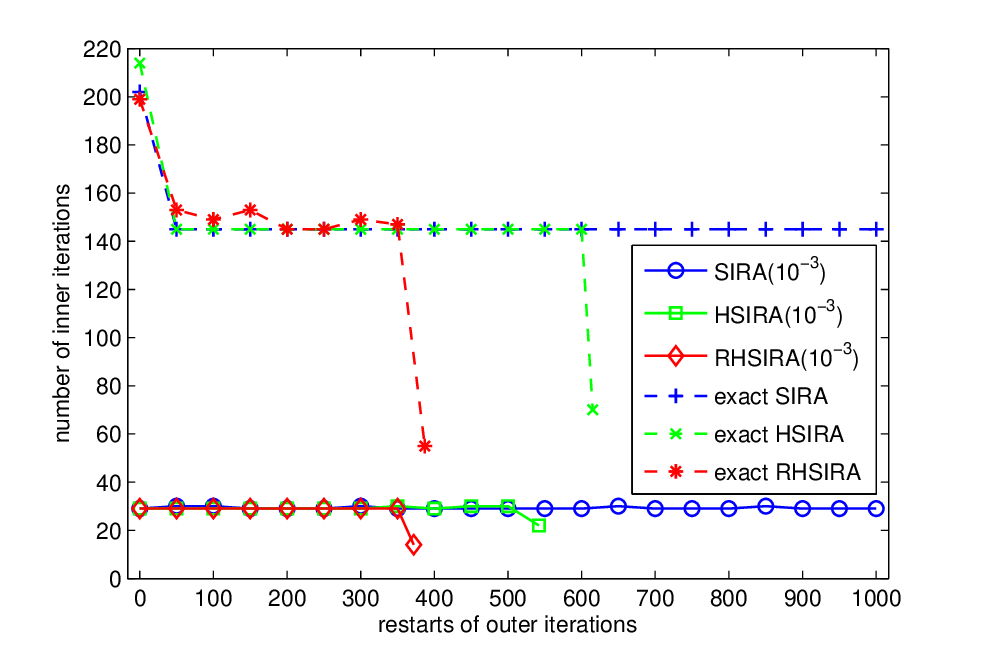}
\end{minipage}
\caption{\emph{Example 4. dw8192 with $\sigma=-60$.}}
\label{fig_dw8192}
\end{figure}

\noindent\textbf{Example 4.} This unsymmetric eigenproblem dw8192 of $n=8192$
arises from dielectric channel waveguide problems \cite{matrixmarket}.
We test this example with $\sigma=-60$ and ${\bf Max}=1000$.
The computed eigenvalue is $\lambda\approx -87.795$.
The preconditioner $M$ is the incomplete LU factorization of
$A-\sigma I$ with drop tolerance $10^{-3}$.
Figure~\ref{fig_dw8192} and Tables~\ref{tab_dw8192} report the results.

\begin{table}[!htb]
\centering
\begin{small}
\begin{tabular}{c||c|r|r|r||c|r|r|r}
\hline
Accuracy &  Algorithm  &  $I_{restart}$  & $I_{inner}$ & $P_{0.1}$ &
Algorithm  &  $I_{restart}$  & $I_{inner}$ & $P_{0.1}$\\
\hline\hline
& SIRA & ${\bf Max}$      & $29262$ & $39\%$
& JD & ${\bf Max}$   & $32222$ & $39\%$   \\
$\tilde{\varepsilon}=10^{-3}$
& HSIRA & $542$   & $15953$ & $76\%$
& HJD & $402$   & $11848$ & $81\%$ \\
& RHSIRA & $372$  & $10813$ & $96\%$
& RHJD & $317$  & $9221$ & $96\%$\\
\hline
& SIRA & ${\bf Max}$   & $38502$ & $12\%$
& JD & ${\bf Max}$   & $39413$ & $12\%$    \\
$\tilde{\varepsilon}=10^{-4}$
& HSIRA & $555$   & $17416$ & $51\%$
& HJD & $497$   & $15691$ & $57\%$    \\
& RHSIRA & $386$  & $11595$ & $93\%$
& RHJD & $429$  & $12884$  & $93\%$  \\
\hline
& SIRA & ${\bf Max}$  & $145435$ & $-$
& JD & ${\bf Max}$  & $151453$ & $-$   \\
exact
& HSIRA & $615$   & $89636$  & $-$
& HJD & $617$   & $90056$ & $-$      \\
& RHSIRA & $387$  & $59170$ & $-$
& RHJD & $383$  & $59166$ & $-$ \\
\hline
\end{tabular}
\end{small}
\caption{\emph{Example 4. dw8192 with $\sigma=-60$.}}
\label{tab_dw8192}
\end{table}

All the algorithms gave results that are typically similar to those corresponding
counterparts obtained for Example 3. The comments and conclusions are very
analogous, and we thus omit them.

\section{Conclusions}\label{sec:conc}

The standard SIRA and JD methods may not work well for the interior eigenvalue
problem. The standard Ritz vectors may converge irregularly, and it is also
hard to select correct Ritz pairs to approximate the desired interior eigenpairs.
So the Ritz vectors may not be good choices for restarting,
causing that restarted algorithms may perform poorly. Meanwhile, the Ritz vectors
appear to expand subspaces poorly. In contrast, the harmonic Ritz vectors are
more regular and more reliable approximations
to the desired eigenvectors, so that they may
expand the subspaces better and generate more accurate a sequence of
subspaces when restarting the methods.
Due to the optimality, the refined harmonic Ritz vectors are
generally more accurate than the harmonic ones and are better choices
for expanding subspaces or restarting the methods.
Most importantly, we have
proved that the harmonic and refined SIRA and JD methods generally mimic
their exact counterparts well provided that all inner linear systems
are solved with only low or modest accuracy. To be practical, we have presented
the restarted harmonic and refined harmonic SIRA and JD algorithms. Meanwhile,
we have designed practical stopping criteria for inner iterations.

Numerical experiments have confirmed our theory. They have indicated that
the restarted harmonic and refined harmonic SIRA and JD algorithms are much
more efficient than the restarted standard SIRA and JD algorithms for the
interior eigenvalue problem. Furthermore,
the refined harmonic algorithms are much better than the harmonic ones.
Importantly, the results have demonstrated
that each inexact algorithm behaves like its exact counterpart
when all inner linear systems are solved with only low or modest accuracy.
In addition, numerical results have confirmed that each exact or inexact
SIRA type algorithm is equally as efficient as the corresponding JD type
algorithm for the same $\tilde\varepsilon$.

Our algorithms are designed to compute only one interior eigenvalue and its
associated eigenvector. This is a special nature of SIRA and JD type methods:
they compute only one eigenpair each time. If more than one eigenpairs
are of interest, we may extend the algorithms to this situation in a number
of ways. For example, we may introduce some deflation techniques
\cite{stewart2001eigensystems,vandervorst2002eigenvalue} into them. Also,
we may change the target $\sigma$ to a new one, to which the second desired
eigenvalue is the closest, and apply the algorithms to find it and
its associated eigenvector. Proceed this way until
all the desired eigenpairs are found. Such kind of algorithms is
under developments.

\begin{small}
\bibliographystyle{siam}

\end{small}

\end{document}